\documentclass[11pt,reqno]{amsart}

\usepackage{amsfonts}
\usepackage{amssymb}           
\usepackage{amsthm}      
\usepackage{amsbsy}      
\usepackage{amscd}      
\usepackage[dvips]{graphicx}
\usepackage{pstricks}
\usepackage{mathrsfs}

\newtheorem{teo}{Theorem}
\newtheorem{lema}[teo]{Lemma}

\theoremstyle{definition}
\newtheorem{definitie}[teo]{Definition}

\newtheorem{remarca}[teo]{Remark}

\newcommand\ze{{\mathbb Z}}

\newcommand\mcg{{\mathcal M}}
\newcommand\hcg{{\mathcal H}}
\newcommand\stab{{\mathcal K}}

\newcommand\braid{{\rm B}}
\newcommand\purebraid{{\rm P}}

\newcommand\fund{{\pi_1}}
\newcommand\rp{{\rm P}}
\newcommand\rd{{\rm D}}

\newcommand\rto{{\longrightarrow}}
\newcommand\dto{{\downarrow}}
\newcommand\eq{{\ =\ }}
\newcommand\vzero{{\overrightarrow{v_0}}}

\begin{document}
\title[Presentation of the handlebody group of genus 2]{A simple presentation of the handlebody group of genus 2$^*$}
\author{Clement Radu Popescu}      
\address{Inst. of Math. "Simion Stoilow", P.O. Box 1-764, RO-014700, 
Bucharest, Romania}%
\email{Radu.Popescu@imar.ro}%

\thanks{$^*$This is part of my Ph.D.Thesis I defended at \textsc{Columbia University} in the spring of 2001}

\begin{abstract}
For genus $g = 2$ I simplify Wajnryb's presentation of the handlebody group.
\end{abstract}

\maketitle

\section{Introduction}
\label{intro}

Let $S_{g,h,k}$ be a closed surface of genus $g$, $h$ boundary components and 
$k$ distinguished points. I'll use the notation $S_g$ for $S_{g,0,0}$. The mapping class group, $\mcg_{g,h,[k]}$ 
is the group of all isotopy classes of orientation preserving 
homeomorphisms which keep the boundary and the distinguished points 
pointwise fixed. $\mcg_{g,h,k}$ is the group of all isotopy classes of orientation preserving homeomorphisms which keep the boundary
pointwise fixed, and is also fixing the set of distinguished points eventualy permuting them. We have
$\mcg_{g,h,[k]}\hookrightarrow \mcg_{g,h,k}$. I'll use the notation $\mcg_g$ for $\mcg_{g,0,0}$. There are some important instances of these groups. For example the braid group $\braid_n$ is shown in
\cite{Birman}, theorem 1.10, to be isomorphic with $\mcg_{0,1,n}$, and so the pure braid group $\purebraid_n$ is isomorphic with $\mcg_{0,1,[n]}$. 

For a long time it was an open problem to obtain a presentation for $\mcg_{g,1,0}$. In \cite{MacCool}  MacCool proved, using purely algebraic methods, that $\mcg_{g,1,0}$ is finitely presented for any genus $g$. Hatcher and Thurston in \cite{Hatcher-Thurston} 
made a crucial breakthrough in the subject developing an algorithm for obtaining a finite presentation for $\mcg_{g,1,0}$. Using this algorithm, Harer in \cite{Harer} obtained a
finite, explicit presentation. Finaly it was Wajnryb, who in \cite{Wajnryb_1} gave a simple presentation for $\mcg_{g,0,0}$ and $\mcg_{g,1,0}$.

Using similar techniques as in \cite{Hatcher-Thurston}, he found in 
\cite{MR2000a:20075} a presentation of the handlebody group. The handlebody group is the subgroup of 
$\mcg_g$ formed by all the isotopy classes 
of orientation preserving homeomorphisms of $S_g$, which extend to 
the entire handlebody $H_g$. I'll denote it by $\hcg_g$. The presentation obtained by Wajnryb in \cite{MR2000a:20075} is long and complicated. 

In this note I will give a simple presentation for $\hcg_2$ starting from Wajnryb's presentation. Such a simplification for higher genera doesn't work. 
 
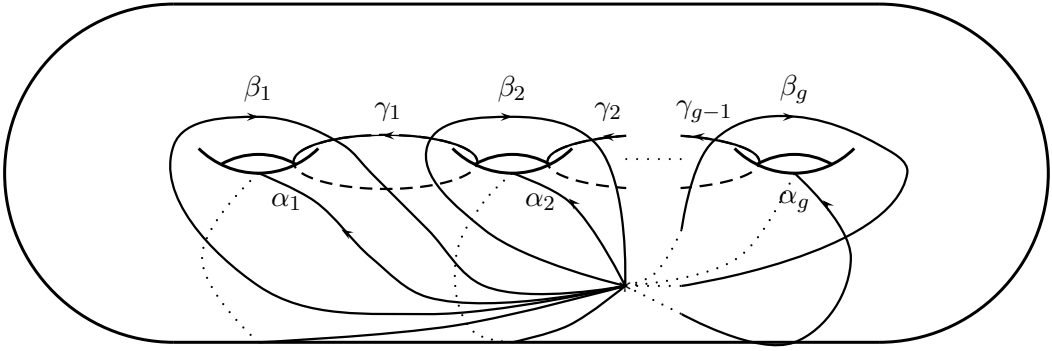
\begin{figure}[h]
\begin{pspicture}(.35,2)(20,7)
\psset{xunit=.75cm,yunit=.75cm}

\rput(10.5,4){$\ast$}


\psline[linewidth=.4mm](2.5,9)(15,9)
\psline[linewidth=.4mm](2.5,3)(15,3)

\psarc[linewidth=.4mm](2.5,6){2.25cm}{90}{270}
\psarc[linewidth=.4mm](15,6){2.25cm}{-90}{90}

\psarc[linewidth=.4mm](4,7.5){1.125cm}{225}{315}
\psarc[linewidth=.4mm](4,5){1cm}{60}{120}
\psarc[linewidth=.4mm](8.5,7.5){1.125cm}{225}{315}
\psarc[linewidth=.4mm](8.5,5){1cm}{60}{120}
\psarc[linewidth=.4mm](13.5,7.5){1.125cm}{225}{315}
\psarc[linewidth=.4mm](13.5,5){1cm}{60}{120}

\pscurve[linewidth=.3mm](10.5,4)(7,3.5)(4,4)(2.5,6.5)(4,7)(5.5,6.5)
(6,6)(7,5)(8,4)(10.5,4)
\psline[linewidth=.3mm]{<-}(4,7)(3.95,7)

\pscurve[linewidth=.3mm](10.5,4)(8,5)(7,6.5)(8.5,7)(10,6.5)(10.5,4)
\psline[linewidth=.3mm]{->}(8.45,7)(8.5,7)

\pscurve[linewidth=.3mm](11.5,5)(12,6.5)(13,7)(15,6.5)(15.5,6)(11.5,4)
\psline[linewidth=.3mm]{->}(13.45,7)(13.5,7)
\pscurve[linewidth=.3mm,linestyle=dotted](10.5,4)(11,4.2)(11.5,5)
\psline[linewidth=.3mm,linestyle=dotted](10.5,4)(11.5,4)



\pscurve[linewidth=.3mm](10.5,4)(7,3.8)(6,4.5)(5,5.5)(4,6)
\pscurve[linewidth=.3mm,linestyle=dotted](4,6)(3,4.5)(4,3)
\pscurve[linewidth=.3mm](4,3)(7.5,3.3)(10.5,4)
\psline[linewidth=.3mm]{->}(5.55,4.95)(5.47,5)


\pscurve[linewidth=.3mm](10.5,4)(9.5,5.5)(8.5,6)
\pscurve[linewidth=.3mm,linestyle=dotted](8.5,6)(8,5.5)(7.5,4.5)(7.8,3.2)(8.5,3)
\pscurve[linewidth=.3mm](8.5,3)(9.5,3.3)(10.5,4)
\psline[linewidth=.3mm]{->}(9.55,5.5)(9.5,5.55)


\pscurve[linewidth=.3mm](11.5,3.5)(13.5,3)(14.5,4.5)(13.5,6)
\pscurve[linewidth=.3mm,linestyle=dotted](13.5,6)(12.5,4.5)(10.5,4)
\psline[linewidth=.3mm,linestyle=dotted](11.5,3.5)(10.5,4)
\psline[linewidth=.3mm]{->}(14.05,5.47)(14,5.52)

\rput(4.5,5.5){$\alpha_1$}
\rput(9,5.5){$\alpha_2$}
\rput(13.5,5.5){$\alpha_g$}
\rput(4,7.5){$\beta_1$}
\rput(8.5,7.5){$\beta_2$}
\rput(13.5,7.5){$\beta_g$}
\rput(6.3,7.1){$\gamma_1$}
\rput(10.2,7.1){$\gamma_2$}
\rput(11.9,7.1){$\gamma_{g-1}$}

\psellipse[linewidth=.3mm,linestyle=dashed](6.25,6.2)(1.65,.5)
\psclip{\psframe[linestyle=none](4.6,6.2)(7.9,6.7)}
\psellipse[linewidth=.3mm](6.25,6.2)(1.65,.5) 
\endpsclip
\psline[linewidth=.3mm]{->}(6.25,6.67)(6.2,6.67)
 
\psellipse[linewidth=.3mm,linestyle=dashed](11,6.2)(1.9,.5)
\psclip{\psframe[linestyle=none](9.1,6.2)(12.9,6.7)}
\psellipse[linewidth=.3mm](11,6.2)(1.9,.5) 
\endpsclip

\psline[linewidth=.3mm]{<-}(11.7,6.66)(11.8,6.65)
\psline[linewidth=.3mm]{->}(10.2,6.65)(10.1,6.66)
\psline[linestyle=dotted](10.5,6.25)(11.5,6.25)
\psframe[linestyle=none,fillstyle=solid,fillcolor=white](10.5,6.5)(11.5,7)
\psframe[linestyle=none,fillstyle=solid,fillcolor=white](10.5,5.5)(11.5,6)
\end{pspicture}
\caption{Surface of genus $g$}
\label{fig:surface-g}
\end{figure}

In Figure \ref{fig:surface-g} it is shown 
a system of curves on the surface $S_g$.
The isotopy classes of the curves 
$\alpha_i$'s will be called meridians, and those of $\beta_i$'s will be called longitudes.

\begin{definitie}
\label{Dehn-twist} A ( positive) Dehn-twist with respect to a simple closed curve $\gamma$, denoted by $T_\gamma$, is a homeomorphism 
of the oriented surface $S_g$, which is supported in a regular 
neighborhood of $\gamma$ and is obtained as follows: one cuts open the surface
along $\gamma$ and rotates one end of it with $360^{\circ}$ to the right and then glue back the surface. This is done in such a way that on the boundary and the complement of the regular neighborhood, $T_\gamma$ is the identity map.
\end{definitie}

I'll call $\gamma$ the support of the Dehn-twist.

The effect of a positive $T_{\gamma}$ on any segment which intercepts the curve $\gamma$ transversally in one point is as follows: cut the segment at the interception point and rotate to the right once around $\gamma$.
The following is a well known result (see \cite{MR2000a:20075}):

\begin{lema}
\label{lema:dehnhomeo}
If h is a homeomorphism of the surface $S_g$, and $T_{\gamma}$
is a Dehn-twist, then $T_{h(\gamma)}\eq hT_{\gamma}h^{-1}.$ 
\end{lema}

\begin{remarca}
\label{rem:conj}
I'll use the following notation $h\ast g\eq hgh^{-1}$. So in this
notation  Lemma \ref{lema:dehnhomeo} can be written 
$T_{h(\gamma)}\eq h\ast T_{\gamma}$
\end{remarca}

Roman letters will be used for a Dehn-twist with a support
denoted by the corresponding Greek letter. For the curves 
in  Figure \ref{fig:surface-g} we have $a_i\eq T_{\alpha_i}$,
$b_i\eq T_{\beta_i}$, $c_i\eq T_{\gamma_i}$. 

\begin{definitie}
A meridinal curve on the surface is one which represents a non-trivial 
homotopy class in $\fund (S_g,\ast)$
and bounds a properly embedded disc in the handlebody $H_g$. 
\end{definitie}

Consider $N_{\alpha}$ the normal subgroup of $\fund (S_g,\ast)$
generated by the homotopy classes of the meridians 
$\alpha_1, \ldots ,\alpha_g$ (see Figure \ref{fig:surface-g}).

\noindent Let $ \#:\mcg_g\rto Aut(\fund (S_g,\ast))$
the homomorphism which takes a homeomorphism 
into the automorphism of the 
fundamental group of the surface, and the image of a 
homeomorphism will be denoted with a $\#$-subscript
$(h\mapsto h_{\#}).$
In \cite{MR28:2530} Griffiths shows that $h\in\hcg_g$ if and only if
$h_{\#}(N_{\alpha})\subset N_{\alpha}$. 
A set of generators for $\hcg_g$ was obtained by Suzuki in 
\cite{MR55:6409}. 

Wajnryb has described in \cite{MR2000a:20075}
an algorithm for getting a presentation of the handlebody group.
This is similar to Hatcher--Thurston's algorithm in \cite{Hatcher-Thurston} and 
to Wanjnryb's algorithm in \cite{Wajnryb_1}. 
For the handlebody $H_g$ of genus $g$ there exists an 
associated 2-dimensional complex $X$, called the 
cut-system complex of the handlebody (different from the one
constructed for a surface $S_g$). The vertices
are cut-systems (collection of g disjoint meridinal curves).
Another cut-system is obtained
if we replace one curve with another meridinal curve disjoint from all the curves
in the cut system. Such vertices are joined by an edge. 
Moreover to any triangle we associate a face, thus we obtain a 2-dimensional simplicial complex, called $X$.

$\hcg_g$ acts transitively on the vertices of $X$, 
and this action can be extended to a simplicial action on $X$. Using this action Wajnryb obtains the presentation for $\hcg_g$. 

In Section \ref{Wajnryb}, I will give some details about Wajnryb's presentation  for any g, and in Section \ref{prezentare g=2} the reduction I obtained for the case $\hcg_2$.  

\section{Wajnryb's algorithm}
\label{Wajnryb}

For the convenience of the reader I will describe in this section
in some detail Wajnryb's algorithm. He proved in \cite{MR2000a:20075}, Theorem 13 that the complex $X$, decribed above 
is connected and simply connected.
The 0-skeleton $X^{(0)}$ has a preferred vertex represented by 
the cut-system formed by the collection of the $g$ meridians, one for each handle.
It is denoted $\vzero\eq$ $<\alpha_1,\alpha_2,\cdots,\alpha_g>.$ 
Wajnryb's algorithm has a few steps.

{\bf Step 1:} Find a presentation of 
the stabilizer of $\vzero$ denoted by $\stab$.
The homeomorphisms in $\stab$ either preserve 
$\vzero$ pointwise or permutes the $\alpha_i$'s or changes
their respective orientations.
A presentation for $\stab$ can be found using the following exact sequences.

\vskip -.1in
\begin{equation}
\begin{array}{ccccccccc}
1 & \rto &(\ze/2\ze)^g &\rto &\pm\Sigma_g &\rto &\Sigma_g &\rto &1
\label{dia:signpe}
\end{array}
\end{equation}

\vskip -.5in
\begin{equation}
\begin{array}{ccccccccc}
\label{dia:stab}
1&\rto&\stab_0&\rto&\stab&\rto&\pm\Sigma_g&\rto&1
\end{array}
\end{equation}

In both sequences \eqref{dia:signpe} and \eqref{dia:stab} 
$\pm\Sigma_g$ is the discrete group of signed
permutations. $\pm\Sigma_g$ is the group of permutations of $\{-g, 1-g,\cdots, 1, 2,\cdots g\}$ such that $\sigma(-i)\eq -\sigma(i)$.
The homomorphism $\pm\Sigma_g\rto\Sigma_g$ is the forgetting sign homomorphism and the sequence \eqref{dia:signpe} splits.  
$\stab_0$ is the subgroup of $\stab$, fixing all the $\alpha_i$'s pointwise.
To find a presentation of $\stab_0$, one needs first a presentation 
for $\mcg_{0,2g,0}$ which can be obtained from the following diagram.

\begin{equation}
\begin{array}{ccccccccc}
\label{dia:stabsphere}
 & & & &1& & & & \\
 & & & &\dto& & & & \\
 & & & &\ze& & & & \\
 & & & &\dto& & & & \\
1&\rto&\ze^{2g}&\rto&\mcg_{0,2g+1,0}&\rto&\mcg_{0,1,[2g]}&\rto&1 \\
 & & & &\dto& & & & \\
1&\rto&\fund(S_{0,2g,1},\ast)&\rto&\mcg_{0,2g,1}&\rto&\mcg_{0,2g,0}&\rto&1 \\
 & & & &\dto& & & & \\
 & & & &1& & & & \\
\end{array}
\end{equation}

Using this presentation of $\mcg_{0,2g,0}$ and the following exact 
sequence

\begin{equation}
\begin{array}{ccccccccc}
\label{dia:stab0}
1&\rto&\ze^g&\rto&\mcg_{0,2g,0}&\rto&\stab_0&\rto&1
\end{array}
\end{equation}

one gets a presentation of $\stab_0$. 
 
{\bf Step 2:} Other cut-systems can be obtained
using "translations". Such a translation is given by the homeomorphism
$r_{i,j}$ which changes $\alpha_j$ into $\gamma_{i,j}$ keeping
all the other $\alpha_i$'s $i\neq j$ fixed (see \cite{MR2000a:20075}).
Wajnryb, using connectedness of $X$, proved that the generators
of $\stab$ together with the elements $r_{i,j}$ generate 
$\hcg_g$. 

{\bf Step 3:} There are finitely many edge
orbits modulo the action of $\hcg_g$. To each edge we 
associate its stabilizer $\stab_{i,j}$. Relations are
coming from the conjugations $r_{i,j}hr_{i,j}^{-1}$ 
and $r_{i,j}^{-1}hr_{i,j}\in\stab$ 
for any $h\in\stab_{i,j}$.
For more details see \cite{MR2000a:20075}. 

Using the above described algorithm he was able to find a presentation of the handlebody group which, unfortunately, is very complicated.

\begin{figure}[h]
\begin{pspicture}(0.25,0)(17,6.3)
\psset{xunit=.75cm,yunit=.75cm}

\psframe[linewidth=.4mm,framearc=.3](18,8)

\pscircle[linewidth=.4mm](4,5.5){.5cm}
\pscircle[linewidth=.4mm](4,2.5){.5cm}
\pscircle[linewidth=.4mm](7,5.5){.5cm}
\pscircle[linewidth=.4mm](7,2.5){.5cm}
\pscircle[linewidth=.4mm](14,5.5){.5cm}
\pscircle[linewidth=.4mm](14,2.5){.5cm}
\pscircle*[linewidth=.4mm](10,5.5){.05cm}
\pscircle*[linewidth=.4mm](10.5,5.5){.05cm}
\pscircle*[linewidth=.4mm](11,5.5){.05cm}
\pscircle*[linewidth=.4mm](10,2.5){.05cm}
\pscircle*[linewidth=.4mm](10.5,2.5){.05cm}
\pscircle*[linewidth=.4mm](11,2.5){.05cm}


\psarc[linewidth=.3mm](4,5.5){.75cm}{0}{180}
\psarc[linewidth=.3mm](4,2.5){.75cm}{180}{360}

\psarc[linewidth=.3mm](4,2.5){.63cm}{90}{270}
\psarc[linewidth=.3mm](7,2.5){.63cm}{-90}{90}

\psarc[linewidth=.3mm](4,5.5){.85cm}{-90}{90}
\psarc[linewidth=.3mm](7,2.5){.75cm}{0}{180}

\psarc[linewidth=.3mm](4,5.5){.95cm}{45}{225}
\psarc[linewidth=.3mm](7,2.5){.9cm}{-135}{45}

\psline[linewidth=.3mm](3,2.5)(3,5.5)
\psline[linewidth=.3mm](5,2.5)(5,5.5)

\psline[linewidth=.3mm](4,3.34)(7,3.34)
\psline[linewidth=.3mm](4,1.66)(7,1.66)

\psline[linewidth=.3mm](3.03,4.7)(6.1,1.7)
\psline[linewidth=.3mm](4.9,6.4)(7.9,3.3)

\psline[linewidth=.3mm,linearc=.3cm](4,4.35)(2.5,4.35)(2.5,1)(6,1)(6,2.5)
\psline[linewidth=.3mm,linearc=.3cm](4,6.63)(1.5,6.63)(1.5,.5)(8,.5)(8,2.5)


\rput(4,5.5){$\partial_{-1}$}
\rput(4,2.5){$\partial_{1}$}
\rput(7,5.5){$\partial_{-2}$}
\rput(7,2.5){$\partial_{2}$}
\rput(14,5.5){$\partial_{-g}$}
\rput(14,2.5){$\partial_{g}$}
\rput(5.6,4.7){$\delta_{-1,1}$}
\rput(.85,4){$\delta_{-1,2}$}
\rput(5.9,3.6){$\delta_{1,2}$}
\rput(7.7,4){$\gamma_1$}
\end{pspicture}
\caption{A disc with $2g$ holes}
\label{fig:holeddisc}
\end{figure}
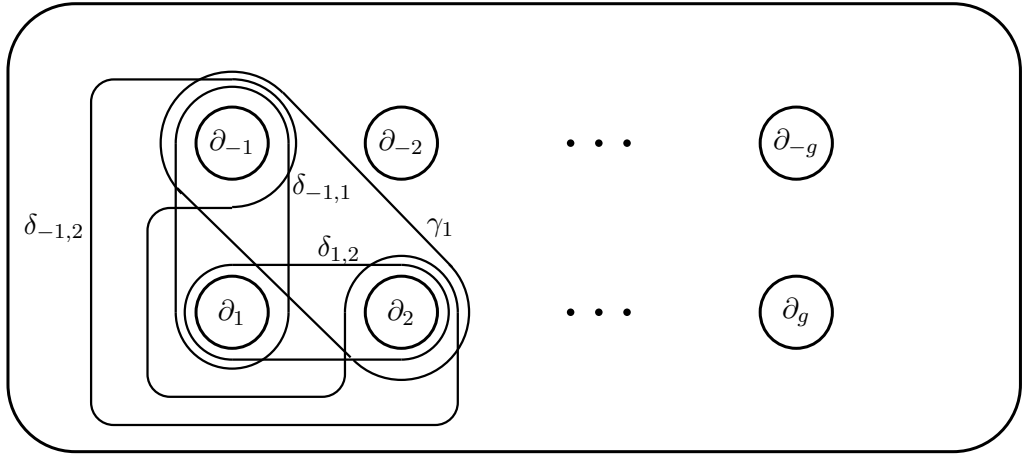

Let me give a brief description of the homeomorphisms 
used in Wajnryb's presentation.  In Figure \ref{fig:holeddisc} 
we see some particular meridinal curves.
The curves denoted by $\gamma_1$  are the same in both Figure 
\ref{fig:surface-g} and \ref{fig:holeddisc}. A positive Dehn-twist is 
considered to be taken in a counterclockwise direction. Looking carefuly at 
Definition \ref{Dehn-twist} one sees that the Dehn-twist  
does not depend on the orientation of the curve.

A half Dehn-twist on $\delta_{-1,1}$ is the twist of the 
first knob and will be denoted by $o$. It changes the orientation of both the
first meridian and longitude. From its definition we see
that $o^2\eq d_{-1,1}$. It is easily seen that 
$o(\delta_{1,2})\eq \gamma_1$. Using Remark \ref{rem:conj}
we get that $o\ast d_{1,2}\eq c_1$. 
Other important homeomorphisms on the surface are
$t_i$ which exchange the meridians $\alpha_i$ and
$\alpha_{i+1}$, fixing the others. Another homeomorphism which exchanges 
$\alpha_i$ and $\alpha_{i+1}$, also exchanging
$\beta_i$ and $\beta_{i+1}$, is $k_i\eq a_ia_{i+1}t_id_{i,i+1}^{-1}$.
I'll mention also $z$, which is
a rotation of the surface about the z-axis 
(the z-axis is considered to pierce the surface in its center
of symmetry as drawn in Figure \ref{fig:surface-g}, 
with positive direction from bellow to above), 
changing the $i$-th hole into the $g-i+1$-st hole, considering that $i > 0$. 
From Figure \ref{fig:holeddisc} it is clear that a curve 
$\delta_{i,j}$ is the one which encloses in one of the two regions determined on the disc, the holes $\partial_i$ and $\partial_j$. 
The Dehn-twists with these supports  can be considered to 
be the generators of the $\purebraid_{2g}\simeq\mcg_{0,1,[2g]}$. 
There is one homeomorphism, $z_j$ which belongs to the stabilizer of an edge of type $(i,j)$ in the 1-skeleton of $X$.
In fact $z_j$ has the form 
$z_j\eq k_{j-1}k_{j-1}\cdots k_{g+j-1}z$
and it is not the conjugation of $z$ by the product of $k_i$'s. 

In the case $g\eq 1$ it is proved by Wajnryb in \cite{MR2000a:20075}, Theorem 14, and detailed in \cite{teza}, Theorem 2.2 that $\hcg_1\simeq\ze\oplus\ze_2$.

In the next statement I will use the same notations as in \cite{MR2000a:20075}, and I will write explicitely Wajnryb's presentation for the case of genus $g\eq 2$.
I will also introduce in the presentation some 
generators together with their defining relations. 

\begin{teo}[Wajnryb] A presentation of $\hcg_2$
is given by:

{\it Generators:} $a_1, a_2, d_{-2-1}, d_{-21}, d_{-22}, 
d_{-11}, d_{-12}, d_{12}, o, o_2, t, r, z, e$.

{\it Defining relations:}
\end{teo}
\label{teo:wgenus2}


$(\rd 1)\hskip .2in o^{-1}t^{-1}o^{-1}\ast d_{12} \eq  d_{-2-1}$ 

$(\rd 2)\hskip .2in t^{-1}o^{-1}\ast d_{12} \eq  d_{-21}$ 

$(\rd 3)\hskip .2in o^{-1}\ast d_{12}  \eq  d_{-12}$

$(\rd 4)\hskip .2in t^{-1}d_{12}\ast d_{-11} \eq d_{-22}$

$(\rd 5)\hskip .2in o_2 \eq  td_{12}^{-1}\ast o $

$(\rd 6)\hskip .2in z \eq a_1^{-1}a_2^{-1}otod_{12}$

$(\rd 7)\hskip .2in e \eq  ozo^{-1}z $

$(\rp 1)\hskip .2in a_1\rightleftarrows a_2\ ;\  
a_i\rightleftarrows d_{kl}$

$(\rp 2.1)\hskip .2in d_{-2-1}^{-1}d_{-21}d_{-2-1} \eq 
d_{-21}d_{-11}d_{-21}d_{-11}^{-1}d_{-21}^{-1}$

$(\rp 2.2)\hskip .2in d_{-2-1}^{-1}d_{-11}d_{-2-1}\eq d_{-21}d_{-11}d_{-21}$

$(\rp 2.3)\hskip .2in d_{-2-1}^{-1}d_{-22}d_{-2-1}\eq 
d_{-22}d_{-12}d_{-22}d_{-12}^{-1}d_{-22}^{-1}$

$(\rp 2.4)\hskip .2in d_{-21}^{-1}d_{-22}d_{-21}\eq
d_{-22}d_{12}d_{-22}d_{12}^{-1}d_{-22}^{-1} $

$(\rp 2.5)\hskip .2in d_{-11}^{-1}d_{-22}d_{-11}\eq d_{-22}$

$(\rp 2.6)\hskip .2in d_{-2-1}^{-1}d_{-12}d_{-2-1}\eq
d_{-22}d_{-12}d_{-22}^{-1}$

$(\rp 2.7)\hskip .2in d_{-21}^{-1}d_{-12}d_{-21}\eq   d_{-22}d_{12}d_{-22}^{-1}d_{12}^{-1}d_{-12}
d_{12}d_{-22}d_{12}^{-1}d_{-22}^{-1}$

$(\rp 2.8)\hskip .2in d_{-11}^{-1}d_{-12}d_{-11} \eq
d_{-12}d_{12}d_{-12}d_{12}^{-1}d_{-12}^{-1}$

$(\rp 2.9)\hskip .2in d_{-2-1}^{-1}d_{12}d_{-2-1}\eq d_{12}$

$(\rp 2.10)\hskip .2in d_{-21}^{-1}d_{12}d_{-21}\eq d_{-22}d_{12}d_{-22}^{-1}$

$(\rp 2.11)\hskip .2in d_{-11}^{-1}d_{12}d_{-11}\eq 
d_{-12}d_{12}d_{-12}^{-1}$

$(\rp 3)\hskip .2in d_{-2-1}d_{-21}d_{-12}d_{-22}
d_{-11}d_{-12}d_{12}\eq a_1^4a_2^4$

$(\rp 4.1)\hskip .2in d_{-11}d_{-12}d_{12}\eq  a_1^2a_2^2$

$(\rp 4.2)\hskip .2in d_{-21}d_{-22}d_{12}\eq a_1^2a_2^2$

$(\rp 4.3)\hskip .2in d_{-2-1}d_{-22}d_{-12}\eq a_1^2a_2^2$

$(\rp 4.4)\hskip .2in d_{-2-1}d_{-21}d_{-11}\eq a_1^2a_2^2$ 

$(\rp 5)\hskip .2in empty\ for\ genus\ 2$

$(\rp 6)\hskip .2in o^2\eq d_{-11}\  ;\  t^2\eq d_{12}d_{-2-1}a_1^{-2}a_2^{-2}$ 

$(\rp 7)\hskip .2in t\ast a_1\eq a_2\ ;\ o\rightleftarrows a_i,\quad i\eq 1,2$

$(\rp 8)\hskip .2in t\rightleftarrows d_{12}\ ;\ otot\eq toto\ ;
\ o\rightleftarrows d_{-22}$

$(\rp 9)\hskip .2in r^2 \eq  a_2^{-4}o_2d_{12}o_2d_{12}^{-1}$

$(\rp 10\ a)\hskip .2in r\ast a_2 \eq d_{12}\ ;\ r\rightleftarrows a_1$ 

$(\rp 10\ b-d)\hskip .2in empty\ for\ genus\ 2 $

$(\rp 10\ e)\hskip .2in r \rightleftarrows e $ 

$(\rp 10\ f)\hskip .2in r\ast d_{12}\eq a_2 $

$(\rp 10\ g)\hskip .2in r\ast d_{-2-1}\eq d_{-11}d_{-12}d_{12}a_1^{-2}a_2^{-1}$

$(\rp 10\ h-k)\hskip .2in empty\ for\ genus\ 2$

$(\rp 11)\hskip .2in rtr \eq  trt $

$(\rp 12)\hskip .2in  empty\ for\ genus\ 2$ 

\vskip .5in

\section{A simple presentation for $\hcg_2$}
\label{prezentare g=2}

In this section I will prove the following:

\begin{teo}
\label{teo:mygenus2}
There is a simple presentation of $\hcg_2$.

Generators: $a_1, a_2, d, o, t, r$.

Relations:
\end{teo}

\begin{center}
$\hskip .2in d \rightleftarrows oto$\ ;\ 
$\hskip .2in odod \eq a_1^2a_2^2$\ ;\ 
$\hskip .2in o^2\rightleftarrows t^{-1}d$\ ;\ 
$\hskip .2in z\eq a_1^{-1}a_2^{-1}otod$\ ;
\end{center}

\begin{center}
$\hskip .2in a_1\rightleftarrows a_2\ ;\ a_i\rightleftarrows d $\ ;\ 
$\hskip .2in t^2\eq d^2a_1^{-2}a_2^{-2}$\ ;\ 
$\hskip .2in o\rightleftarrows a_i\ ;\ t\ast a_1\eq a_2$\ ;\ 
\end{center}

\begin{center}
$\hskip .2in t\rightleftarrows d\ ;\ otot\eq toto $\ ;\ 
$\hskip .2in r^2\eq  d^{-2}a_1^2a_2^{-2}$\ ;\ 
$\hskip .2in r\ast a_2 \eq d_{12}\ ;\ r\rightleftarrows a_1$\ ;\ 
\end{center}

\begin{center}
$\hskip .2in r\rightleftarrows  ozo^{-1}z$\ ;\ 
$\hskip .2in rtr\eq trt$.
\end{center}

\begin{proof}
Use (\rp 4.4) and (\rp 4.1) in (\rp 3) and also the commuting relations (\rp 1) to get 

\[(\rp 3)'\hskip .2in d_{-11}\eq d_{-22} \] 

Replace $d_{-22}$ with $d_{-11}$ in relations (\rp 4.1)--(\rp 4.4), denoted 
(\rp 4.1)'--(\rp 4.4)'. Using the commuting relations (\rp 1) and (\rp4.3)' 
one gets 
\vskip .1in
$(\rp 4.1)''\hskip .2in d_{12}\eq d_{-2-1}$.
\vskip .1in

Rewrite again (\rp 4.2)'--(\rp 4.4)', and use (\rp 4.1)'' to obtain 
(\rp 4.2)''-- (\rp 4.4)''. Cojugate (\rp 4.4)'' with $d_{12}^{-1}$ and modulo the commuting relations (\rp 1) one gets again (\rp 4.2)''.
So (\rp 4.4)'' is redundant. So far (\rp 3)--(\rp 4.4) look like:

\vskip .1in
$(\rp 3)'\quad \hskip .2in d_{-11}\eq d_{-22}$

$(\rp 4.1)''\hskip .2in d_{-2-1}\eq d_{12}$

$(\rp 4.2)''\hskip .2in d_{-21}\eq a_1^2a_2^2d_{12}^{-1}d_{-11}^{-1}$

$(\rp 4.3)''\hskip .2in d_{-12}\eq a_1^2a_2^2d_{-11}^{-1}d_{12}^{-1}$

$(\rp 4.4)''\hskip .2in redundant$
\vskip .1in
Using the above relations in (\rp 2.1)--(\rp 2.11) (the pure braid relations), these all become trivial modulo the commuting relations (\rp 1). This is the main reduction in the presentation, which does not take place for any higher genus. So, for $g\eq 2$ the pure braid relations are redundant modulo (\rp 1) and concequences of (\rp 3)--(\rp 4.4). 

Using Tietze operations, we replace in the remaining relations the expressions for $d_{-2-1}$, $d_{-21}$, $d_{-22}$, $d_{-12}$ and remove these generators together with (\rp 3)''--(\rp 4.4)''. I will remove also generators $d_{-11}$, $o_2$, $e$
together with first part of (\rp 6), (\rd 5) and (\rd 7) replacing first their expressions everywhere else.

At this moment the presentation looks like this:

\vskip .1in
$Generators: a_1, a_2, d_{12}, o, t, r, z$

$Relations:$

$(\rd 1)\hskip .2in o^{-1}t^{-1}o^{-1}\ast d_{12}\eq d_{12}$ 

$(\rd 2)\hskip .2in t^{-1}o^{-1}\ast d_{12}\eq a_1^2a_2^2d_{12}^{-1}o^{-2}$ 

$(\rd 3)\hskip .2in o^{-1}\ast d_{12}\eq a_1^2a_2^2o^{-2}d_{12}^{-1}$

$(\rd 4)'\hskip .2in t^{-1}d_{12}\ast o^2\eq o^2$

$(\rd 6)'\hskip .2in z\eq a_1^{-1}a_2^{-1}otod_{12}$

$(\rp 1)\hskip .2in a_1\rightleftarrows a_2\ ;
\ a_i\rightleftarrows d_{12}\ ;\ 
a_i\rightleftarrows o^2$

$(\rp 6)'\hskip .2in t^2\eq d_{12}^2a_1^{-2}a_2^{-2}$

$(\rp 7)\hskip .2in t\ast a_1\eq a_2\ ;\ o\rightleftarrows a_i$

$(\rp 8)\hskip .2in t\rightleftarrows d_{12}\ ;\ otot\eq toto\ ;\ o\rightleftarrows o^2$

$(\rp 9)'\hskip .2in r^2\eq a_2^{-4}(td_{12}^{-1}\ast o)d_{12}(td_{12}^{-1}\ast o)d_{12}^{-1}$

$(\rp 10a)\hskip .2in r\ast a_2 \eq d_{12}\ ;\ r\rightleftarrows a_1$

$(\rp 10e)'\hskip .2in r\rightleftarrows ozo^{-1}z$

$(\rp 10f)'\hskip .2in r\ast d_{12} \eq a_2$ 

$(\rp 10g)'\hskip .2in r\ast d_{12} \eq 
o^2a_1^2a_2^2o^{-2}d_{12}^{-1}d_{12}a_1^{-2}a_2^{-1}$ 

$(\rp 11)\hskip .2in rtr\eq trt$
\vskip .1in

Let $d_{12}\eq d$ ( it is the only $d_{ij}$ left). In the above (\rd 1)
is equivalent with $(\rd 1)' \hskip .2in oto\rightleftarrows d$.

Rewrite (\rd 2) and use (\rd 1)' to get (\rd 2)':

\begin{center}
$
t^{-1}o^{-1}dot\eq a_1^2a_2^2d^{-1}o^{-2}
\Leftrightarrow t^{-1}o^{-1}doto\eq a_1^2a_2^2d^{-1}o^{-1}
\Leftrightarrow 
$
\end{center}
\begin{center}
$
(using\ \rp 1)\Leftrightarrow t^{-1}o^{-1}otod\eq a_1^2a_2^2d^{-1}o^{-1}\Leftrightarrow  
od\eq a_1^2a_2^2d^{-1}o^{-1}\Leftrightarrow
(\rd 2)'\hskip .2in odod\eq a_1^2a_2^2.
$
\end{center}

Conjugating the above with $d$ and using (\rp 1) we get 
a redundant relation:

\begin{center}
$(\rd 2)''\hskip .2in dodo\eq a_1^2a_2^2\eq odod$.
\end{center}

Rewrite (\rd 3) using (\rd 2)'' as follows:
\begin{center}
$
o^{-1}do\eq a_1^2a_2^2o^{-2}d^{-1}\Leftrightarrow 
o^{-1}dodo\eq a_1^2a_2^2o^{-1}
(using\ (\rd 2)'')\Leftrightarrow
o^{-1}a_1^2a_2^2\eq a_1^2a_2^2o^{-1} (trivial\ because\ of\ (\rp 7)).
$
\end{center}

Relation(\rd 4)' is equivalent with $t^{-1}d_{12}\rightleftarrows o^2$.

In (\rp 1) the last part, $o^2\rightleftarrows a_i$, is redundant from the
second part of (\rp 7) which is $o\rightleftarrows a_i$. 

In (\rp 8) we can get rid of $o\rightleftarrows o^2$. 

Lastly in (\rp 10g)', because of the commuting relations (\rp 1) and (\rp 7) we'll get (\rp 10f)'. So (\rp 10g)' is redundant. 

Rewrite (\rp 9)'. Use for this (\rp 8), (\rd 2)'',
(\rp 1) to get $r^2\eq d^{-2}a_1^{2}a_2^{-2}$.

Because (\rp 10f)' and the above one easily gets that the 
first part of (\rp 10a) is redundant. 

This is the presentation in the statement of Theorem \ref{teo:mygenus2}

\end{proof}

\end{document}